\documentclass[11pt]{amsart}

\usepackage[english]{babel}
\usepackage{amssymb}
\usepackage{fourier}
\usepackage{mathrsfs}
\usepackage{enumerate}
\usepackage{color}
\usepackage{ifpdf}
\usepackage{tikz}
\usepackage{tikz-cd}
\usetikzlibrary{decorations.pathmorphing,arrows, intersections}
\usepackage[initials]{amsrefs}
\usepackage{enumitem}

\usepackage{caption}
\usepackage{subcaption}
\usepackage{comment}

\newcommand{\bbR}{{\mathbf R}}

\newcommand{\bbZ}{{\mathbf Z}}

\newcommand{\GL}{\operatorname{GL}}

\newcommand{\SL}{\operatorname{SL}}
\newcommand{\PSL}{\operatorname{PSL}}

\newcommand{\Homeo}{\operatorname{Homeo}}

\newcommand{\overto}[1]{{\buildrel{#1}\over\longrightarrow}}

\newcommand{\setdef}[2]{ \left\{ {#1}\ /\ {#2} \right\} }

\newtheorem{theorem}{Theorem}[section]

\newtheorem{mthm}{Theorem}

\newtheorem{lemma}[theorem]{Lemma}
\newtheorem{claim}[theorem]{Claim}

\newtheorem{proposition}[theorem]{Proposition}

\newtheorem{definition}[theorem]{Definition}
\newtheorem{defn}[theorem]{Definition}

\newtheorem{example}[theorem]{Example}

\newtheorem{remark}[theorem]{Remark}

\numberwithin{equation}{section}

\begin{document}
\title[Marked length pattern rigidity]{Marked length pattern rigidity\\ for arithmetic manifolds}

\author{Yanlong Hao}
\address{University of Illinois at Chicago, Chicago, Illinois, USA}
\email{yhao23@uic.edu}
 
\date{\today}
    
\maketitle

\begin{abstract}
    In this paper, we prove a cocycle version of marked length spectrum rigidity. There are two consequences.  The first is marked length pattern rigidity for arithmetic hyperbolic locally symmetric manifolds.  The second is strengthen marked length spectrum rigidity for surfaces and closed locally symmetric manifolds.
\end{abstract}

\section{Introduction}
Let $(M,g)$ be a closed Riemannian manifold whose sectional 
curvatures are all negative, 
$\Gamma=\pi_1(M)$ its fundamental group, and $C_\Gamma$ 
the set of conjugacy classes 
of non-trivial elements of $\Gamma$.
Then classes in $C_\Gamma$ correspond to free non-trivial homotopy 
classes of loops in $M$,
and each such class contains a unique loop of minimal $g$-length 
- a closed geodesic.
\emph{Marked length spectrum} is the function
\[
    \ell_g:C_\Gamma\overto{} [0,\infty)
\]
that assigning the $g$-length $\ell_g(\langle \gamma\rangle)$
of the closed geodesic corresponding to the conjugacy
class $\langle \gamma\rangle\in C_\Gamma$ of
$\gamma\in\Gamma\setminus\{1\}$.

\medskip

\textit{Marked Length Spectrum Rigidity} Conjecture (cf. Burns and Katok
\cite{burns1985manifolds}) states that 
function $\ell_g$ determines $(M,g)$, up to isometry. 

\medskip

This conjecture was proved for surfaces by Otal \cite{otal1990spectre} and, 
independently but in greater generality, 
by Croke \cite{croke1990rigidity} slightly later. 
For higher dimensions, Katok \cite{katok1988four} had a short proof for metrics 
in fixed conformal class, 
the proof is given in dimension 2 only, but can be easily extended 
to all dimensions. 
Beside this, for
higher dimension Hamenstädt \cite{hamenstadt1997cocycles} 
proved the conjecture in the case where $(M,g)$ 
is a locally symmetric space by applying  the
famous entropy rigidity work of Besson-Courtois-Gallot \cite{besson1995entropies}. 
For general negatively curved metrics the problem is largely open. 
Guillarmou and Lefeuvre showed in \cite{guillarmou2019marked} that the marked length spectrum of a Riemannian manifold $(M,g)$ with Anosov geodesic flow and non-positive curvature 
locally determines the metric $g$.  

\medskip

In this paper we consider a slightly different problem, that concerns rigidity of what
we call \textbf{marked length pattern}. 
Given, say closed, negatively
curved Riemannian manifold $(M,g)$ with fundamental group $\Gamma=\pi_1(M)$ 
and conjugacy classes $C_\Gamma$, 
the \textit{length pattern} is the equivalence relation
$R_g$ on $C_\Gamma$ given by equality of $g$-length:
\[
    R_g=\setdef{(c_1,c_2)\in C_\Gamma\times C_\Gamma}{\ell_g(c_1)=\ell_g(c_2)}.
\]
We emphasize that the relation $R_g$ describes equality of $g$-lengths, but not their values. 
\begin{mthm}\label{b}
    Let $(M,g_0)$ be a closed arithmetic locally symmetric manifold of rank 1, and let $g$ 
    be an arbitrary negatively curved Riemannian metric on $M$.
    Then $R_{g_0}\subset R_g$ only if $(M,g)$ is isometric to $(M,\lambda g_0)$ 
    for some $\lambda>0$ by an isometry isotopic to identity. 
\end{mthm}

\medskip

On the other hand. there are strengthen marked length spectrum rigidity.

A.~Katok \cite{katok1988four} proved that marked length spectrum
in a fixed homology class determines the negatively curved metric 
on the surface
in a fixed conformal class. In \cite{gogolev2020abelian}, 
Gogolev and Rodriguez Hertz showed that the marked length spectrum restricted
to the set of conjugacy classes represented by homologically 
trivial geodesics is
enough to uniquely determine the marked length spectrum. 
Noelle \cite{S} showed the same is true for compliment of 
a `small' set when the manifold is a closed surface.
The following Theorem extends some of these results,
including \cite{gogolev2020abelian}. 

\begin{mthm}\label{d}
Let $(M,g_1)$ and $(M,g_2)$ be two arbitrary closed negatively curved
Riemannian metrics on a manifold $M$ with fundamental group $\Gamma$.
Let $H$ be a subgroup of $\Gamma$ such that the limit set of $H$ is all of $\partial \widetilde{M}$. 

Then, $\ell_{g_1}=\ell_{g_2}$ on classes from $H$ 
only if $\ell_{g_1}=\ell_{g_2}$ on all of $\Gamma$. 
Moreover, if $M$ is a surface or $(M, g_1)$ is a locally 
symmetric rank one space, then $(M,g_1)$ is isometric to $(M,g_2)$.
\end{mthm}

For example, this applies to any non-trivial $H$ which is
a normal subgroup $\Gamma$, or a normal subgroup of normal subgroup, etc.

In a forthcoming paper, joint work with Alexander Furman, 
we give a generalization of Noelle's work \cite{S} using 
a different method.  

\medskip

Theorem \ref{b} and \ref{d} deal with different perspectives. 
However, the key of their proofs relies on the same framework:
marked length spectrum rigidity for $B$-cocycles. To state it,  
let us start with the following

\begin{defn}\label{D:hyp-element}
Let $X$ be a compact space. An element $\phi\neq Id_X$ in $\Homeo(X)$ 
is called \textbf{hyperbolic} if $\phi$ has two distinct
fixed points $\phi^-,\phi^+\in X$ and $\phi^n$, $n\to+\infty$, 
contract $X\setminus\{\phi^-\}$ towards $\phi^+$ uniformly on 
compact subsets, and $\phi^{n}$, $n\to-\infty$,
contracts $X\setminus\{\phi^+\}$ towards $\phi^-$.
\end{defn}

This is often called \textit{north-south dynamics}.
Note that the points $\phi^{\pm}\in X$ 
are uniquely determined by  a hyperbolic $\phi$.
Note also that if $\phi$ is hyperbolic, then so is $\phi^{-1}$
and its attracting/repelling points are $(\phi^{-1})^\pm=\phi^\mp$.

\begin{defn}\label{D:geom-bnd}
Let $\Gamma$ be a topological group. A nontrivial compact 
Hausdorff $\Gamma$-space $X$ is called a \textbf{geometric boundary} if
\begin{enumerate}
    \item $\Gamma$ acts on $X$ minimally,
    \item Every $\gamma\in\Gamma\setminus\{e\}$ is hyperbolic,
    \item There are $\Gamma$-quasi-invariant measures 
    $\mu$, $\mu'$ on $X$ such that the $\Gamma$-action on 
    $(X\times X,\mu\times \mu')$ is ergodic. 
\end{enumerate}
If in addition, $\mu=\mu'$, we call $X$ a 
symmetric geometric boundary. 
\end{defn}


There are many examples of symmetric geometric boundaries. A typical example is a hyperbolic group $G$ and its boundary $\partial G$, or any non-elementary subgroup of a hyperbolic group 
with its limit set (using Kaimanovich \cite{Kaim} and Bader-Furman \cite{bader2014boundaries}), 
or some subgroups of acylindrically hyperbolic groups and their limit sets (using
Maher-Tiozzo \cite{maher2018random}).


\medskip

In this paper, we mainly consider the applications to the marked length spectrum related problems. 

Let $\Gamma=\pi_1(M)$ be the fundamental group of a closed negatively
curved Riemannian manifold $M$. The fundamental group $\Gamma$ acts
continuously on $\partial \widetilde{M}$. In fact, $\Gamma$-space
$\partial \widetilde{M}$ is a geometric boundary. Fixed a base point,
there is Busemann cocycle 
$\alpha:\Gamma\times \partial \widetilde{M}\rightarrow \bbR$.
The cocycle $\alpha$ encodes many geometric properties of the metric.
In particular, for any nontrivial $\gamma\in \Gamma$,
$\alpha(\gamma,\gamma^+)=\ell_g(\langle \gamma\rangle)$. 

This inspires us to define the \emph{marked length spectrum function} 
for general cocycle $\beta$ on a geometric boundary $X$ 
by setting
\begin{equation}
    \ell_\beta(\langle \gamma\rangle)=\beta(\gamma,\gamma^+)
\end{equation}
for all $\gamma\in \Gamma\setminus\{1\}$.

It is more straightforward to define the marked length spectrum 
of $\beta$ by a pair of numbers 
$\ell_\beta(\langle \gamma\rangle)=(\beta(\gamma,\gamma^+),\beta(\gamma,\gamma^-))$.
However, our definition is enough for the special type of cocycles,
namely \textit{$B$-cocycles} that we define below, where we shall see $\beta(\gamma,\gamma^-)=-\beta(\gamma,\gamma^+)$
and $\ell_\beta(\langle \gamma^{-1}\rangle)=\ell_\beta(\langle \gamma\rangle)$
(see section \ref{MLS} below).

\begin{defn}\label{D:B-cocycles}
Let $X$ be a $\Gamma$-space.
A cocycle $c:\Gamma\times X\rightarrow \bbR$ is a  \textbf{B-cocycle} if
\begin{enumerate}
    \item The map $x\mapsto c(\gamma,x)$ 
    is continuous for all $\gamma\in\Gamma$.
    \item There exists a continuous function 
    $C:X\times X\setminus\Delta\rightarrow \bbR$ such that $$C(\gamma x,\gamma y)-C(x,y)=c(\gamma,x)+c(\gamma,y)$$ for all $\gamma\in \Gamma$, $x\neq y\in X$ where $\Delta$ is the diagonal $\setdef{(x,x)}{x\in X}$ of $X\times X$.
\end{enumerate}
We also call $(c,C)$ a $B$-cocycle when it is necessary to point out the function $C$. 
\end{defn}

Denote by $Z^1_c(\Gamma,X,\bbR)$ the vector space of all
continuous cocycles $c:\Gamma\times X\to\bbR$,
and by $B^1_c(\Gamma,X,\bbR)$ the subspace 
consisting of cocycles
of the form $c(\gamma,x)=\varphi(\gamma.x)-\varphi(x)$
for some continuous function $\varphi:X\to\bbR$.
Denote by $H^1_c(\Gamma,X,\bbR)=Z^1_c/B^1_c$ 
the associated cohomology.
Equivalently, we can view this as cohomology $H^1(\Gamma,C(X,\bbR))$
with coefficients in $C(X,\bbR)$ viewed as a $\Gamma$-module.

Note that the collection of all $B$-cocycles over $\Gamma$-space $X$
forms a vectors subspace of $Z_c^1(\Gamma,X,\bbR)$
that contains $B_c^1(\Gamma,X,\bbR)$, hence one can talk about
cohomology classes of $B$-cocycles, or just $B$-\textbf{classes}.
$B$-classes form the kernel of the map 
\[
    [i]:H^1_c(\Gamma,X,\bbR)\to H^1_c(\Gamma,X\times X\setminus\Delta,\bbR)
\]
induced by $i\alpha(g,(x,y))=\alpha(g,x)+\alpha(g,y)$.
We have the following rigidity results.

\begin{mthm}\label{a}
Let $\Gamma$-space $X$ be a geometric boundary, 
and $\alpha, \beta: \Gamma\times X\rightarrow \bbR$
two $B$-cocycles. Then
$\ell_\alpha=\ell_\beta$ if and only if $[\alpha]=[\beta]$
in $H_c^1(\Gamma,X,\bbR)$.
\end{mthm}

\begin{remark}\label{16}
Even though lattices in semisimple Lie groups $G$ and their Furstenberg boundaries $X=G/P$ are not necessarily geometric boundaries, there is a similar definition of B-cocycles in this case where the diagonal $\Delta$ in (2) is replaced by 
the set of pairs of points not in general position. 
B-cocycle can be defined also for the general acylindrical groups and their limit sets.

Theorem \ref{a} is still true in both cases. However, the proof is not cover in our framework. We will prove it for general acylindrical groups in the appendix. The proof of higher rank lattices need a little bit more work, but follows form the same construction. The detailed proof for higher rank lattices will be given in other place. 
\end{remark}

 In the setting of negatively curved space $(M,g)$ 
 and $\Gamma=\pi_1(M)$ acting on $X=\partial\widetilde M$ the 
 theorem implies that the cross-ratios on the boundaries is determined by the Busemann cocycles when restricted to Busemann cocycles coming from various Riemannian metrics. Hence it is a weaker statement than the marked length spectrum rigidity and is well-known in this case. In fact, for the Busemann cocycle of a closed negatively curved manifold, there is a corresponding length cocycle for its geodesic flow. Fix one visual metric on the boundary associated to a specified negatively curved metric. All pullbacks of Busemann cocycles for different metrics are H\"{o}lder. And in this case, Theorem \ref{a} can be deduced from a Theorem of Liv{\v{s}}ic \cite{livvsic1972cohomology}.  

Even though it is not new in many cases, Theorem
\ref{a} provides us a new viewpoint to see more patterns behind the marked length spectrum of an arithmetic locally symmetric space. It is not the length themselves, but rather some identities between them determine the metric up to homothecy. Theorem \ref{b} follows from this observation. Indeed, Theorem of Liv{\v{s}}ic \cite{livvsic1972cohomology} about H\"{o}lder cocycles and our construction together is enough to prove Theorem \ref{b}.

\medskip

It is natural to ask: what is the case for non-arithmetic lattices? The short answer is: they do not have marked length pattern rigidity in general.  

\begin{mthm}\label{f}
For a general hyperbolic surface $(S,g)$, the relation $R_g$ is the minimal relation among all relations from hyperbolic metrics. In other words, $R_g$ is a sub-relation of $R_{g'}$ for any hyperbolic metric $g'$ on $S$.
\end{mthm}

Here, hyperbolic metric means a Riemaniann metric with constant curvature $-1$. Hyperbolic surface is a surface with a finite volume complete hyperbolic metric.  

We do not know any example of higher dimension locally symmetric manifold which do not have the marked length pattern rigidity. By Mostow rigidity, the locally symmetric structure is determined by the fundamental group itself. We need more tools to discuss general metrics.

\medskip
Similar to Remark \ref{16}, Theorem \ref{b} and \ref{d} have weaker generalization to finite volume manifolds and orbifolds. Since all these need a slightly different setup and the results are weaker, we discuss them in the appendix.  
\medskip

The paper is organized as follows. In Section \ref{c} we recall some basic facts about Gromov-hyperbolic spaces, Patterson-Sullivan measure and Busemann cocycles. We also define the $B$-cocycles there. In Section \ref{4}, we give some examples of geometric boundaries and show some very basic properties. In Section \ref{MLS}, we prove Theorem \ref{a}. Section \ref{6} contains a proof of Theorem \ref{d}. Section \ref{extension} is concerned on extension of $B$-cocycles from arithmetic lattice $\Gamma<G$ to the Lie group  $G$. In Section \ref{l}, we apply the results in Section 6 to negatively curved Riemaniann manifolds, this yields Theorem \ref{b}. In section \ref{hypernolic surfaces}, we discuss the marked length spectrum of hyperbolic surfaces. And  give examples of hyperbolic metric without marked length patter rigidity. The Appendix \ref{appendix} is devoted to show marked length spectrum rigidity for cocycles of general acylindrical groups. 

\medskip

I would like to acknowledge and give my thanks to my supervisor Alexander Furman who made this work possible. His guidance and advises help me work over all the challenges.  

\section{Gromov-Hyperbolic spaces and Patterson-Sullivan measures}\label{c}

We use \cite{bridson2013metric} as a general reference for Gromov hyperbolic spaces, and here only add some specifics that are not covered there. 

Let $X$ be a proper CAT(-1) space. Fix a base-point $p\in X$. The Gromov product $(\cdot,\cdot)_p$ on $X\times X$ extends continuously and therefore canonically to a function $\overline{X}\times \overline{X}\rightarrow [0,\infty]$, and  for each $\epsilon \in(0,1]$ 
the expression $d_\epsilon(\xi,\eta)=e^{-\epsilon (\xi,\eta)_p}$ defines a metric on the boundary $\partial X=\overline{X}\setminus X$ 
compatible with the topology of $\overline{X}$ (cf. Bourdon \cite{bourdon2016cohomologie}).
We call them visual metrics. In particular, $\xi=\eta$ on $\partial X$ if and only if $(\xi, \eta)_p=\infty$.

Let $(M,g)$ be a strictly negatively curved complete Riemannian manifold, the universal cover of $M$ is a Gromov-hyperbolic space, in fact, a CAT$(-k)$ space, for some $k>0$. Its visual boundary is the same as its Gromov boundary. By renormalizing the metric, we see that the Gromov product extends continuously to the boundary. There are still visual metrics $d_\epsilon$ for all $\epsilon\leq \epsilon_0$.     


\subsection{Patterson-Sullivan measures}\label{e}
\

Let $\Gamma$ be a non-elementary discrete group of isometries of a connected, simply-connected, complete Riemannian manifold $(M,g)$ with the curvature $\kappa\leq -a^2<0$. 
Denote by $\overline{M}=M\cup \partial M$ the compactification of $M$, it is homeomorphic to a closed ball.

In the case of dimension 2 and of constant curvature, Patterson (\cite{patterson1976limit}) gave a construction of a family of measures on $\partial M$ indexed by points $x\in M$. It has been extended to higher dimensions by Sullivan (\cite{sullivan1979density}). Generalizing their work, a number of authors, including Coornaert, Albuquerque, and Knieper constructed the so-called Patterson-Sullivan measures $\nu_x$ on $\partial M$ in variable curvature case. See \cite{albuquerque1999patterson} and \cite{knieper1997asymptotic} for the construction of Patterson-Sullivan measures.

\medskip

Let us discuss some key properties of Patterson-Sullivan measures $\nu_x$. 

\begin{enumerate}
    \item Equivariance: $\nu_{\gamma x}=\gamma_* \nu_x$ for all $\gamma \in \Gamma$.
    \item  Explicit Radon-Nikodym derivative: For all $x,y\in M$,
         \begin{equation}\label{E2}
            \frac{d\nu_x}{d\nu_y}(\xi)=e^{h(g)B_{x,y}(\xi)},
        \end{equation}
    where $B_{x,y}(\xi)$ is the Busemann function on $M$ and $h(g)$ in the volume entropy of $M$. 

For points $x,y\in M$ and $\xi\in \partial M$, the function 
$B:M\times M\times \partial M\rightarrow \bbR$ is defined by
\[
    B_{x,y}(\xi)=\lim_{t\rightarrow \infty}(d_X(y,\gamma_\xi(t))-t)
\]
where $\gamma_\xi$ is the unique geodesic ray with $\gamma(0)=x$ and $\gamma(\infty)=\xi$. 

\begin{remark}
 In general, $\nu_x$ are just finite measures, not probability measures. 
 \end{remark}
\begin{remark}
The Radon-Nikodym derivative is defined almost everywhere with respect to $\nu_y$. 
However, using the fact that $B$ is continuous, and the support of $\nu_y$ is the limit set of $\Gamma$, 
we can assume Equation (\ref{E2}) is true for all points in the limit set of $\Gamma$.
\end{remark}
\end{enumerate}
\begin{enumerate}[start=3]
    \item\label{(3)} Marked length spectrum is determined by $\nu_x$: When $M/\Gamma$ is a manifold, $\Gamma$ is torsion free. Let $\gamma\in \Gamma$ be a hyperbolic element. There is a unique attracting fixed point $\gamma^{+}$ of $\gamma$ on $\partial M$. 

\begin{claim} 
    $B_{x,\gamma x}(\gamma^{+})=\ell_g(\langle \gamma\rangle)$ for all $x\in M$.
\end{claim}

\begin{proof}
By equation \ref{E2}, we have 
\[
    B_{x,\gamma x}(\gamma^+)=B_{y,\gamma y}(\gamma^+)+B_{x,y}(\gamma^+)
    -B_{\gamma x,\gamma y}(\gamma^+).
\]

Since $\gamma$ is an isometry, by the definition of Busemann function, 
\[
    B_{\gamma x,\gamma y}(\gamma^+)=B_{\gamma x,\gamma y}(\gamma\gamma^+)=B_{x,y}(\gamma^+).
\]
Hence $B_{x,\gamma x}(\gamma^+)$ is independent of the choice of $x$. Let $x$ be a point on the axis of $\gamma$. 
The claim follows. 
\end{proof}

For a parabolic element $\eta$, there is only one fixed point $\xi$ on $\partial M$, 
and 
\[
    B_{x,\eta x}(\xi)=0.
\]
Hence we set $\ell_g(\langle \eta\rangle)=0$.
\item\label{(4)} Bowen-Margulis-Sullivan geodesic current on $\partial M\times \partial M$: Fix a base point $p\in M$. 
Let $\epsilon>0$ be small enough so that $d_\epsilon$ is a metrics for all $x\in M$. 
A direct computation shows that 
\[
    d\mu(\xi,\eta)=d_\epsilon(\xi,\eta)^{-2h(g)/\epsilon}\cdot dv_p(\xi)dv_p(\eta)
\]
defines a $\Gamma$-invariant measure on $\partial M\times \partial M$. In fact, $d\mu$ is the Bowen-Margulis-Sullivan geodesic current. 
It is independent of the choice of $p$ and $\epsilon$.
\end{enumerate}
\subsection{Cocycles}
\

Recall that when $X$ is a $\Gamma$-space, a continuous function $\alpha:\Gamma\times X\rightarrow \bbR$ is called a \emph{cocycle} if for all $\gamma$, $\eta\in \Gamma$, and $x\in X$, 
$$\alpha(\gamma\eta,x)=\alpha(\gamma,\eta x)+\alpha(\eta,x).$$

For any map $k: X\rightarrow \bbR$, define the \emph{coboundary} of $k$ by 
\[
    d k(\gamma, x)=k(\gamma x)-k(x)\qquad(\gamma\in \Gamma,\ x\in X).
\]
All coboundaries are cocycles.  
\begin{definition}
Two cocycle $\alpha$, $\beta$ are called equivalent, if there is a continuous map $k:X\rightarrow \bbR$ such that $$\alpha-\beta=d k.$$
\end{definition}
It is not surprising that we can define cohomology group of a $\Gamma$-space by identify all equivalent cocycles as a cohomological class. Let $Z$ be the set of cocycles, and $B$ the set of coboundaries of continuous maps.
$Z$ is an abelian group, and $B$ is a subgroup of $Z$. The first cohomology group $H^1_c(\Gamma,X,\bbR)$ is just $Z/B$. In fact, it is the same thing as the usual group cohomology of $\Gamma$ with coefficient $C(X,\bbR)$, all continuous map from $X$ to $\bbR$ viewed as a $\Gamma$-module. 

Notice that there is a more general definition of cocycles. The target could be any second countable group, instead of $\bbR$, and the cohomology is not necessarily a group in this setting. There are also Borel cocycles and measurable cocycles when the $\Gamma$-space $X$ admit a non-singular measure, see \cite{zimmer2013ergodic}. Recall the following definition in \cite{zimmer2013ergodic}.
\begin{definition}
Two Borel cocycles $\alpha$, $\beta$ are called strictly equivalent if there is a Borel function $\phi$ such that 
$$\alpha-\beta=d \phi.$$
\end{definition}

It looks a little bit strange to define strictly equivalence in this way. It is the notion taken in \cite{zimmer2013ergodic}. The notion of equivalence of cocycles here is different form that in \cite{zimmer2013ergodic}. We will work for continuous map in this paper in most cases. Only refer to Borel cocycles and Borel maps for some technical issues.

\medskip

Let $G_0$ be a closed subgroup of $G$. Then $G$ acts on $G/G_0$ via left translation. Let $H$ be a second countable group. The following Proposition was showed in \cite{zimmer2013ergodic}.
\begin{proposition}
There is a bijection 
\[
    \left\{Borel\  Cocycles \ G\times G/G_0\rightarrow H\right\}\rightarrow {\rm{Hom}}(G_0, H)
\]
between strict equivalence classes of Borel cocycles and conjugacy classes of homomorphisms.
\end{proposition}

\subsection{$B$-Cocycles}
\

There are Busemann-cocycles closely related to the Patterson-Sullivan measures. Fix a base point $o\in M$, the Busemann-cocycle $B(\gamma, \xi)$ is given by
$$B(\gamma, \xi)=B_{\gamma^{-1}o,o}(\xi)=\frac{1}{h(g)}\ln \frac{d(\gamma^{-1}_\ast \nu_o)}{d\nu_o}(\xi),$$
for all $\gamma\in \Gamma$ and $\xi\in \partial M$.

By (\ref{(3)}) in \ref{e}, the marked length spectrum is determined by the Busemann-cocycle. 
And (\ref{(4)}) in \ref{e} is equivalent to
$$-2(\gamma\xi,\gamma\eta)_o+2(\xi,\eta)_o=B(\gamma,\xi)+B(\gamma,\eta)$$ for all $\gamma\in \Gamma$, $\xi,\eta\in \partial M$ and $\xi\neq \eta$.

In general, for any $\Gamma$-space $X$.
\begin{defn}
A cocycle $c:\Gamma\times X\rightarrow \bbR$ is a  \textbf{$B$-cocycle} if
\begin{enumerate}
    \item $c(\gamma,\xi)$ is continuous for all $\gamma\in\Gamma$.
    \item There exist continuous function $C:X\times X\setminus\Delta\rightarrow \bbR$ 
    such that $$C(\gamma x,\gamma y)-C(x,y)=c(\gamma,x)+c(\gamma,y)$$ for all 
    $\gamma\in \Gamma$, $x\neq y\in X$ where $\Delta$ is the diagonal of $X\times X$.
\end{enumerate}
We say $(c,C)$ is a $B$-cocycle when it is necessary to point out the function $C$. 
\end{defn}

For any $\Gamma$-space $X$, there is diagonal action of $\Gamma$ on $X\times X$. 
Consider the map 
\[
    [i]: H^1_c(\Gamma, X,\bbR)\rightarrow H^1_c(\Gamma, X\times X\setminus\Delta,\bbR),
\]
defined by $i(\alpha)(\gamma,(x,y))=\alpha(\gamma,x)+\alpha(\gamma,y)$.
Then $B$-cocycles are exactly the cocycles represent classes in $\text{Ker}([i])$.
\begin{lemma}\label{R:C-determines-c}
    In a B-cocycle $(c,C)$, the function $C:X\times X\setminus \Delta\to \bbR$ determines the cocycle 
    $c:\Gamma\times X\overto{} \bbR$.
\end{lemma}
\begin{proof}
    Let $h(x,y,z)=C(x,y)+C(x,z)-C(y,z)$ for three pairwise different points $x$, $y$, $z\in X$. 
    Then 
    \[
        2c(\gamma, x)=h(\gamma x,\gamma y,\gamma z)-h(x,y,z).
    \]
\end{proof}
%

\section{Geometric Boundaries}\label{4}
In this section, we give some basic properties of geometric boundaries which will be useful in next section.

Recall the definition of geometric boundaries.
\begin{defn}
Let $\Gamma$ be a topological group. A nontrivial compact 
Hausdorff $\Gamma$-space $X$ is called a \textbf{geometric boundary} if
\begin{enumerate}
    \item $\Gamma$ acts on $X$ minimally,
    \item Every $\gamma\in\Gamma\setminus\{e\}$ is hyperbolic,
    \item There are $\Gamma$-quasi-invariant measures 
    $\mu$, $\mu'$ on $X$ such that the $\Gamma$-action on 
    $(X\times X,\mu\times \mu')$ is ergodic. 
\end{enumerate}
If in addition, $\mu=\mu'$, we call $X$ a 
symmetric geometric boundary. 
\end{defn}

Note that the support of $\mu$ is a $\Gamma$-invariant closed subset of $X$. Hence $\mu$ has full support since the $\Gamma$-action is minimal. Same is true for $\mu'$.

In Lemma \ref{R:C-determines-c}, it was showed that for a $B$-cocycle $(c, C)$, the function $C$ determines the cocycle $c$. On the other hand, we have the following:
\begin{lemma}Let $\Gamma$-space $X$ be a geometric boundary, and $(c,C)$ a $B$-cocycles on $X$. Then the cocycle $c$ determines $C$ up to a constant.
\end{lemma}
\begin{proof}Assume $(c,C)$ and $(c,C')$ are both $B$-cocycles. 

By definition $C-C'$ is a $\Gamma$-invariant function. 
Since $\Gamma$ is $\nu\times \nu'$-ergodic, $\nu\times \nu'$ has full support and $C-C'$ is continuous, 
$C-C'$ is a constant function. 
\end{proof}

We give some examples of geometric boundaries.

\begin{example}\label{E:BM-boundary}
Let $\Gamma$ be the fundamental group of a closed negatively curved Riemannian manifold $(M,g)$, and $\partial \widetilde{M}$ the boundary of the Riemanian universal cover $\widetilde{M}$.  The $\Gamma$-space $\partial \widetilde{M}$ is a symmetric geometric boundary.
\end{example}

\begin{enumerate}
  \item Since the limit set of $\Gamma$ coincide with $\partial \widetilde{M}$, the action is minimal.

  \item It is well-known and can be deduced by the classification of isometries of hyperbolic spaces.  

  \item Consider the Patterson-Sullivan measures $\nu_x$ on $\partial \widetilde{M}$. There is the Bowen-Margulis-Sullivan geodesic current on $\partial^2 \widetilde{M}$ which is in the measure class $[\nu_x\times \nu_x]$. It is well-known that the Bowen-Margulis-Sullivan geodesic current is $\Gamma$-ergodic. 
\end{enumerate}

\begin{example}
Now, with same setting as in Example~\ref{E:BM-boundary}. Let $H\lhd\Gamma$ be a normal subgroup such that 
$\Gamma/H$ is a finite extension of $\bbZ$ or $\bbZ^2$. 
Then $\partial \widetilde{M}$ is a symmetric geometric $H$-boundary.
\end{example}
\begin{enumerate}
    \item  Since $H$ is normal, the limit set of $H$ is the same as the limit set of $\Gamma$. The action is minimal.

\item $H$ is not a trivial group. 

\item The Bowen-Margulis-Sullivan geodesic current is a $H$-invariant ergodic measure (cf. Guivarc'h \cite{MR959369}).
\end{enumerate}

Now we give some basic properties of geometric boundaries. All statements are well known. We give some short proofs for completeness.

\begin{lemma}\label{Q}
Geometric boundaries have infinite many points.
\end{lemma}
\begin{proof}
Let $X$ be a geometric boundary of $\Gamma$. There exists a hyperbolic element $\gamma$ by (2) in Definition \ref{D:geom-bnd}. Since $\gamma$ fixes $\gamma^+$, $\gamma^-$, and $\gamma$ acts on $X$ non-trivially, $X$ contains at least 3 points. For any $x\neq \gamma^+$, 
$x\neq \gamma^-$, $$\lim_{n\rightarrow +\infty}\gamma^n x=\gamma^+.$$
It is clear that $\gamma^n x\neq \gamma^+$ for all $n$. This completes the proof
\end{proof}

It is not hard to show that in fact, geometric boundaries are perfect spaces. 

\begin{lemma}\label{P}
Let $\Gamma$-space $X$ be a geometric boundary. Let $\gamma$ be a hyperbolic element. There exist $\theta\in \Gamma$ such that $\theta \gamma^+\neq \gamma^-$, $\theta \gamma^-\neq \gamma^-$.
\end{lemma}
\begin{proof}
Assume not. Take any element $\theta$ in $\Gamma$. $\theta \gamma^+=\gamma^-$ or $\theta \gamma^-=\gamma^-$.

If $\theta \gamma^+=\gamma^-$, $\theta^2 \gamma^+=\theta\gamma^-\neq\gamma^-$. It follows that $\theta^2\gamma^-=\gamma^-$. Hence
$\theta \gamma^-=\theta^{-1}\gamma^-=\gamma^+$.

The $\Gamma$-orbit of $\gamma^-$ is just the set of two points $\{\gamma^+, \gamma^-\}$. The action of $\Gamma$ is not minimal by Lemma \ref{Q}. This is a contradiction.
\end{proof}
A direct computation shows that
\begin{lemma}\label{R}
Conjugations of a hyperbolic element are hyperbolic. And $$(\theta\gamma\theta^{-1})^\pm=\theta\gamma^\pm$$
when $\gamma$ is hyperbolic.
\end{lemma}
\begin{lemma}\label{S}
Let $\gamma, \eta\in \Gamma$ with $\gamma$ hyperbolic, then $\eta\gamma^n$ is hyperbolic for all but possibly one $n$. Furthermore, if $\eta \gamma^+\neq \gamma^-$, then 
$$\lim_{n\rightarrow +\infty}(\eta\gamma^n)^+=\eta\gamma^+,$$
$$\lim_{n\rightarrow +\infty}(\eta\gamma^n)^-=\gamma^-.$$
\end{lemma}
\begin{proof}
First claim follows by the fact that $\gamma$ acts on $X$ non-trivially and has an infinite order. Hence at most one of $\eta\gamma^n$ acts trivially on $X$. 

Let $U$, $V$ be two open sets such that $U\cap V=\emptyset$, $\eta \gamma^+\in U$, $\gamma^-\in V$.

By definition, there exists $M$ such that for all $n\geq M$, $\gamma^n(X-V)\subset \eta^{-1}U$. 
It implies $\eta\gamma^n(X-V)\subset U$.  Therefore $\eta\gamma^n$ is hyperbolic for $n$ big enough. 
And for the hyperbolic element $\eta\gamma^n$, when $x\neq (\eta\gamma^n)^-$,
\[
    \lim_{k\rightarrow +\infty}(\eta\gamma^n)^k x=(\eta\gamma^n)^+.
\]
Hence $(\eta\gamma^n)^+\in U$ for $n$ big enough. Since $U$ can be chosen arbitrary small, 
\[
    \lim_{n\rightarrow \infty}(\eta\gamma^n)^+=\eta\gamma^+.
\]
Similarly, for $n$ big enough, $\gamma^{-n}\eta^{-1}(X-U)\subset V$ and $(\gamma^{-n}\eta^{-1}))^+\in V$.
$$\lim_{n\rightarrow \infty}(\eta\gamma^n)^-=\lim_{n\rightarrow \infty}(\gamma^{-n}\eta^{-1}))^+=\lim_{n\rightarrow \infty}\eta(\eta^{-1}\gamma^{-n}))^+=\gamma^-,$$
since $\eta^{-1}\gamma^{-n}=\eta^{-1}(\gamma^{-n}\eta^{-1})\eta$.
\end{proof}
\section{Marked length spectrum rigidity for $B$-cocycles}\label{MLS}

We prove Theorem \ref{a} in this section. 

Before the proof, we first make the following observation.
Let $(c,C)$ be a $B$-cocycle over a geometric $\Gamma$-boundary $X$. Since 
\[
    c(\gamma,\gamma^+)+c(\gamma,\gamma^-)=C(\gamma \gamma^+,\gamma \gamma^{-})-C(\gamma^+,\gamma^-)=0
\]
for all $\gamma\in \Gamma\setminus \{1\}$.
We have 
\[
    c(\gamma,\gamma^-)=-c(\gamma,\gamma^+).
\]
By 
\[
0=c(1,\gamma^-)=c(\gamma^{-1}\gamma,\gamma^-)=c(\gamma^{-1},\gamma\gamma^-)+c(\gamma,\gamma^-)=c(\gamma^{-1},\gamma^-)+c(\gamma,\gamma^-),
\]
\[\ell_c(\langle\gamma^{-1}\rangle)=c(\gamma^{-1},\gamma^-)=c(\gamma,\gamma^+)=\ell_c(\langle \gamma\rangle).
\]

We restate Theorem \ref{a} here for reader's convenience.
\begin{theorem}
Let $\Gamma$-space $X$ be a geometric boundary and $\alpha, \beta: \Gamma\times X\rightarrow \bbR$
two $B$-cocycles. Then $\ell_\alpha=\ell_\beta$ if and only if  $[\alpha]=[\beta]$ in $H^1_c(\Gamma,X,\mathbb{R})$.
\end{theorem}
\begin{proof}
The `if part' is trivial.

For `only if part', let $\delta=\alpha-\beta$. It is again a $B$-cocycle and $\ell_\delta=0$. We need to show $\delta=d\varphi$ for some continuous function $\varphi$.
Let us assume the pair $(\delta,f)$ is a $B$-cocycle.

The proof has 3 steps.

\medskip

\textbf{Step 1:} The cocycle $\delta$ is bounded.

Choose any hyperbolic element $\gamma\in \Gamma$ with fixed points $\gamma^+$ and $\gamma^{-}$. Let $\eta\in \Gamma$ be any element such that $\eta \gamma^+\neq \gamma^-$. We have
\[\delta(\eta \gamma^n,\gamma^+)=\delta(\eta,\gamma^n\gamma^+)+\delta(\gamma^n,\gamma^+)=\delta(\eta,\gamma^+)\]
for all $n$ by assumption.

Now for $n$ big enough, 
\[
\delta(\eta,\gamma^+)=\delta(\eta \gamma^n, \gamma^+)=\delta(\eta \gamma^n, \gamma^+)+\delta(\eta \gamma^n, (\eta\gamma^n)^-).
\]
Here, the last equality follows from the fact $\delta(\eta \gamma^n, (\eta\gamma^n)^-)=-\delta(\eta \gamma^n, (\eta\gamma^n)^+)=0$. Since $(\delta,f)$ is a $B$-cocycle,
\[
\delta(\eta \gamma^n, \gamma^+)+\delta(\eta \gamma^n, (\eta\gamma^n)^-)=f(\eta \gamma^n \gamma^+, \eta \gamma^n (\eta\gamma^n)^-))-f(\gamma^+,(\eta\gamma^n)^-)).
\]
The right hand side is just $f(\eta\gamma^+, (\eta\gamma^n)^-))-f(\gamma^+,(\eta\gamma^n)^-))$.

Taking limit as $n\rightarrow +\infty$, by Lemma \ref{S},
\[\delta(\eta,\gamma^+)=f(\eta\gamma^+,\gamma^-)-f(\gamma^+,\gamma^-).\]

Choose $\theta\in \Gamma$ such that $\theta(\gamma^+)\neq \gamma^-$, $\theta(\gamma^-)\neq\gamma^-$. This is possible by Lemma \ref{P}. The same argument for $\theta \gamma \theta^{-1}$ instead of $\gamma$ as before implies when $\eta\theta \gamma^+\neq \theta \gamma^-$, 
\[
\delta(\eta,\theta\gamma^+)=f(\eta\theta\gamma^+,\theta\gamma^-)-f(\theta\gamma^+,\theta\gamma^-).
\]

Fix two open neighbourhood $U_1$ and $U_2$ of $\gamma^-$ and $\theta\gamma^-$, respectively such that $U_1\cap U_2=\emptyset$, $\gamma^+\notin U_1$, $\theta\gamma^+\notin U_1$, $\theta \gamma^+\notin U_2$. The existence of $U_1$ and $U_2$ is guaranteed since $X$ is Hausdorff. 

Recall that $f$ is continuous, in particular, $f(\cdot,\gamma^-):X-\{\gamma^-\}\rightarrow \bbR$ is continuous. Since $X$ is compact, and $U_1$ is open, $X-U_1$ is compact. The extreme value theorem provides $M_1$ with $|f(x,\gamma^-)|\leq M_1$ for all $x\notin U_1$. Similarly, there exists $M_2$ with $|f(y,\theta\gamma^-)|\leq M_2$ for all $y\notin U_2$. 

Set $M=\max\{M_1, M_2\}$. Now we are ready to show that $\delta$ is bounded.

First, we show $\delta$ has an uniform bound on a point. For any element $\zeta\in \Gamma$.

Case 1: $\zeta \gamma^+\notin U_1$. Then 
\[
|\delta(\zeta,\gamma^+)|=|f(\zeta\gamma^+,\gamma^-)-f(\gamma^+,\gamma^-)|\leq 2M.
\]

Case 2: $\zeta \gamma^+\in U_1$. Then 
$\zeta \gamma^+\notin U_2$, and
\[
|\delta(\zeta\theta^{-1},\theta\gamma^+)|=|f(\zeta\gamma^+,\theta\gamma^-)-f(\theta\gamma^+,\theta\gamma^-)|\leq 2M,
\]
\[
|\delta(\theta,\gamma^+)|=|f(\theta\gamma^+,\gamma^-)-f(\gamma^+,\gamma^-)|\leq 2M.
\]

Hence
\[
|\delta(\zeta,\gamma^+)|=|\delta(\zeta\theta^{-1}\theta,\gamma^+)|=|\delta(\zeta\theta^{-1},\theta \gamma^+)+\delta(\theta,\gamma^+)|\leq 4M.
\]

In both cases, there is an uniformly bounds: $|\delta(\zeta,\gamma^+)|\leq 4M$ for all $\zeta\in \Gamma$.

Second, $\delta$ is bounded on a $\Gamma$-orbit. By cocycle identity, 
\[
\delta(\rho,\varrho \gamma^+)=\delta(\rho\varrho,\gamma^+)-\delta(\varrho,\gamma^+)
\]
is bounded by $8M$ for all $\rho,\varrho\in \Gamma$. 

Since the orbit of $\gamma^+$ is dense and $\delta$ is a continuous cocycle, $\delta$ is globally bounded. 

\medskip

\textbf{Step 2:} Find $\phi$ as a Borel function.

It is well known that a bounded cocycle is a coboundary $d\varphi'$ for some Borel function $\varphi'$, see, for example, \cite{furman2002coarse}. It is trivial to show that $f(x,y)-\varphi'(x)-\varphi'(y)$ is a $\Gamma$-invariant Borel function.

\medskip

\textbf{Step 3:} $\phi$ is essentially continuous.

 $\Gamma$ acts ergodicly on $(X\times X, \mu\times \mu')$. We know $f(x,y)-\varphi'(x)-\varphi'(y)$ is constant $\mu\times \mu'$-a.e.. A Fubini type argument shows that there is a $\mu'$-conull subset $A$, such that for all $y_0\in A$, $f(x,y_0)-\varphi'(x)-\varphi'(y_0)$ is constant $\mu$-a.e.. $\varphi'(x)$ is $\mu$-a.e. equal to the continuous function $\varphi(x):=f(x,y_0)-\varphi'(y_0)$ up to a constant.  Technically, $f(y_0,y_0)$ is undefined. Hence $\varphi$ is just continuous on $X\setminus\{y_0\}$. However, run the same argument for $y_1\neq y_0$ in $A$. And by gluing the two functions together, $\varphi$ is continuous on $X$. 
 
 Since $d\varphi$ and $\delta$ are two continuous cocycles $\mu$-a.e. identical, they are the same cocycle, and $\delta=d\varphi$.
\end{proof}
 
 The marked length spectrum $\ell$ defines a map from $H^1(\Gamma,X,\bbR)$ to $C(\Gamma)^\bbR$. What we have showed so far is that when restricted on $\text{Ker}([i])$, this map is injective. The image of this map is mysterious to us. $\ell$ is not injective in general even for geometric boundaries.

\section{Proof of Theorem \ref{d}}\label{6}
\begin{proof}
Let $(M,g)$ be a closed Riemannian manifold of negative curvature, $\Gamma=\pi_1(M)$ its fundamental group,
$X=\partial\widetilde{M}$ -- boundary of the universal covering, and $H<\Gamma$ a subgroup that acts minimally on $X$.

Let $\mu$ be a symmetric, generating measure on $H$ that has a finite first moment with respect to the distance
$d_{\widetilde{g}}$ on $\widetilde{M}$, i.e. $\sum_{h\in H}\mu(h)\cdot d_{\widetilde{g}}(e,h)<\infty$.
It follows from the work of Kaimanovich \cite{Kaim}*{\S 7.3} that the Poisson boundary of $(H,\mu)$
is realized on $X$, more precisely that there is a (unique) $\mu$-stationary probability measure $\nu$ on $X$,
so that $(X,\nu)$ is the Poisson boundary for $(H,\mu)$.

It follows from Bader and Furman \cite{bader2014boundaries}*{Theorem 2.7} 
that $\nu\times\nu$ on $X\times X$ is $H$-ergodic. 
Therefore $X$ is a symmetric $H$-boundary.

\medskip







Let $g_1$ and $g_2$ be two negatively curved Riemannian metrics on $M$,
and let $\beta_1,\beta_2\in Z^1_c(\Gamma,X,\bbR)$ be the Busemann cocycles associated
with the lifted metrics $\widetilde{g}_1$ and $\widetilde{g}_2$ on $\widetilde{M}$.
Recall that these are $B$-cocycles: there exist (geometrically defined) continuous functions
$f_1,f_2:X\times X\setminus\Delta\to \bbR$ so that 
\[
    \beta_i(\gamma,x)+\beta_i(\gamma,y)=f_i(\gamma x,\gamma y)-f_i(x,y)\qquad (\gamma\in\Gamma,\ x\ne y\in X,\ i=1,2).
\]
The restrictions $\bar\beta_i: H\times X\to\bbR$ are still $B$-cocycles, namely $(\bar\beta_i,f_i)$.
By assumption, $\ell_{\bar\beta_1}=\ell_{\bar\beta_2}$. 
By marked length spectrum rigidity of cocycles, there exist continuous map $\varphi$ such that $\bar\beta_1-\bar\beta_2=d \varphi$. 

Direct calculation show that the continuous function $X\times X\setminus\Delta\to\bbR$ given by
\[
    f_1(x,y)-f_2(x,y)-\varphi(x)-\varphi(y)
\]
is $H$-invariant.
Since $H$ is $\nu\times\nu$-ergodic, while $\nu$ has full support the function is $\nu\times\nu$-a.e. 
constant. Being continuous, and since $\nu\times\nu$ has full support on $X\times X$, the functions
are actually constant. 
This implies that $\beta_1-\beta_2=d \varphi$.
Hence $\ell_{\beta_1}=\ell_{\beta_2}$ on $\Gamma$. 
In other words, $(M,g_1)$ and $(M,g_2)$ have the same marked length spectrum. 

The results for surface and locally symmetric spaces follow from the marked length spectrum rigidity.
\end{proof}

\section{Extension of $B$-cocycles}\label{extension}
In this section, $\Gamma$ will be a torsion-free uniform arithmetic lattice in a rank 1 simple center-free Lie group $G$.  Let $X$ be the Furstenberg boundary of $G$, which is the same as the Gromov boundary of $G/K$. All cocycles in this section will be understood defined over $X$. Let $G'<G$ be a subgroup. We call a cocycle $\omega: G'\times X\rightarrow \bbR$ a cocycle of $G'$. 

Denote by $\alpha$ the Busemann-cocycle associated to the locally symmetric space $\Gamma\backslash G/K$ 
with base point $[e]$. Let $(\beta,f)$ be another $B$-cocycle of $\Gamma$. 
Define equivalence relation $R_\alpha$ on $C_\Gamma$ by  
\[
    (\langle\gamma_1\rangle,\langle\gamma_2\rangle)\in  R_\alpha \qquad\Longleftrightarrow \qquad \ell_\alpha (\langle\gamma_1\rangle)=\ell_\alpha(\langle\gamma_2\rangle).
\]
We define $R_\beta$ similarly.

The main goal for this section is the following Proposition.
\begin{proposition}\label{Pro:extension}
If $R_\alpha$ is a sub-relation of $R_\beta$, then
there is a cocycle $\bar\beta:G\times X\to\bbR$ which extends $\beta$ to the Lie group $G$.
\end{proposition}

\medskip

There are two steps for this extension. First, we extend $\beta$ to the commensurability subgroup
$\text{Comm}_G(\Gamma)$ of $\Gamma$ in $G$. Then we extend $\beta$ from $\text{Comm}_G(\Gamma)$ to $G$ 
using a density type argument. 
We shall use the following observation. 
Since $\alpha$ is a restriction of a $G$-cocycle $G\times X\to\bbR$, if $\gamma_1,\gamma_2\in \Gamma$ are conjugate in $G$, then $(\langle\gamma_1\rangle,\langle\gamma_2\rangle)\in R_\alpha$
and so $(\langle\gamma_1\rangle,\langle\gamma_2\rangle)\in R_\beta$.

\subsection{Extension to the Commensurability subgroup}
\

Let $s$ be an element of $\text{Comm}_G(\Gamma)$. 
By definition, there are finite index subgroups $\Gamma'$, $\Gamma''$ of $\Gamma$, 
so that $s\Gamma's^{-1}=\Gamma''$. 
Denote by $\beta'$ and $\beta''$ the restrictions
of $\beta:\Gamma\times X\to\bbR$ to $\Gamma'$ and $\Gamma''$ respectively. 
There is another cocycle $s_\ast\beta''$ of $\Gamma'$ defined by 
\[
    s_\ast\beta''(\gamma,\xi)=\beta''(s\gamma s^{-1}, s\xi)\qquad (\gamma\in\Gamma').
\]
Similar construction gives $\alpha'$ and $s_\ast \alpha''$.

Since $\gamma$ and $s\gamma s^{-1}$ are conjugate in $G$, it follows that $(\langle\gamma\rangle,\langle s\gamma s^{-1}\rangle)\in R_{\alpha'}\subset R_{\beta'}$ for $\gamma\in\Gamma'$, and so $\ell_{\beta'}=\ell_{s_\ast\beta''}$ on $\Gamma'$.
Note that $(s_\ast \beta'', f_s)$ is a $B$-cocycle, where $f_s(x,y)=f(s x,s y)$. 
Applying the marked length spectrum rigidity for cocycles of $\Gamma'$, 
there exists a continuous function $\varphi_s$ such that $s_\ast\beta''-\beta'=\delta\varphi_s$.

Now $f(s\xi,s\eta)-f(\xi,\eta)-\varphi_s(\xi)-\varphi_s(\eta)$ is a $\Gamma'$-invariant
continuous function. 
Let $m$ be a $G$-invariant measure on $X\times X\setminus\Delta$. 
Recall that $\Gamma'$ acts on $(X\times X\setminus\Delta,\mu)$ ergodicly by Howe-Moore Ergodicity Theorem, $f(s\xi,s\eta)-f(\xi,\eta)-\varphi_s(\xi)-\varphi_s(\eta)$ is $m$-a.e. a constant. 
Since it is a continuous function and $m$ has full support, the function is constant.  
Let
\[
    \hat\beta(s,\xi)=\varphi_s(\xi)+\frac{f(s\xi,s\eta)-f(\xi,\eta)-\varphi_s(\xi)-\varphi_s(\eta)}{2}.
\]
It is clear that 
\[
    f(s\xi,s\eta)-f(\xi,\eta)=\hat\beta(s,\xi)+\hat\beta(s,\eta).
\]
We now shows that $\hat\beta$ defines a cocycle of $\text{Comm}_G(\Gamma)$. 
Define a new function $h(\xi,\eta,\omega)=f(\xi,\eta)+f(\xi,\omega)-f(\eta,\omega)$ on triples of $X$. 
It is straightforward that for all $s\in \text{Comm}_G(\Gamma)$, $\xi,\eta,\omega\in X$ pairwise different,
\[
    h(s \xi,s \eta,s \omega)-h(\xi,\eta,\omega)=2\hat\beta(s,\xi).
\]
The left hand side is a cocycle. Hence the right hand side is a cocycle, too.
However, the right hand side depends on the first factor only. 

It implies that $\hat\beta(s,\xi)$ is a $B$-cocycle of $\text{Comm}_G(\Gamma)$.

\subsection{Extension to $G$}
\

By the fact that $\Gamma$ is arithmetic, $\text{Comm}_G(\Gamma)$ is dense in $G$ with the Hausdorff topology, 
see \cite{zimmer2013ergodic}.

We will finish the argument by the following lemma.
\begin{lemma}
Let $L<G$ be a dense subgroup, with a B-cocycle $(\hat\beta, f)$ of $L$, then $\hat\beta$ 
extends to a B-cocycle $(\bar\beta,f)$ of $G$.
\end{lemma}
 \begin{proof}
Same as before, define $h(\xi,\eta,\omega)=f(\xi,\eta)+f(\xi,\omega)-f(\eta,\omega)$. 
It is clear that for all $l\in L$, $\xi,\eta,\omega\in X$,
\[
    h(l \xi,l \eta,l \omega)-h(\xi,\eta,\omega)=2\hat\beta(l,\xi).
\]
Fix an arbitrary element $g\in G$, and a sequence $\{l_n\}_{n=1}^{\infty}$ in $L$ so that 
\[
    \lim_{n\rightarrow \infty} l_n=g.
\]
For any element $l\in L$, we can identify $\hat\beta(l,\xi)$ as a function on triples. With this in mind, 
\[
    2\hat\beta(l_n,\xi)-2\hat\beta(l_m,\xi)=h(l_n \xi,l_n \eta, l_n \omega)-h(l_m \xi,l_m \eta, l_m \omega).
\]
The right hand side is the same as
\[
    \Omega(n,m):=[h(l_n \xi,l_n \eta, l_n \omega)-h(g\xi,g\eta,g\omega)]-[h(l_m \xi,l_m \eta, l_m \omega)-h(g\xi,g\eta,g\omega)].
\]
Fix $\xi$, we can choose $\eta$, $\omega$ so that the three points are different pairwise. 
$G$ acts on $X$ continuously. Hence
\[
    \lim_{n,m\rightarrow \infty}\Omega(m,n)=0
\]
uniformly  on compact subset of pairwise different triples. That is, $\hat\beta(l_n,x)$ converges uniformly on compact set. It implies that $\hat\beta(l_n,\xi)$ is a Cauchy sequence that uniformly converge to some continuous function, denote by $\bar\beta(g,\xi)$. 
We have 
\[
    h(g\xi,g\eta,g\omega)-h(\xi,\eta,\omega)=2\bar\beta(g,\xi).
\]
Hence as before $\bar\beta$ is a cocycle, and $\bar\beta(g,\xi)+\bar\beta(g,\eta)=f(g\xi,g\eta)-f(\xi,\eta)$ by definition of $h$.
\end{proof}

\section{Marked length pattern rigidity}\label{l}
We prove marked length pattern rigidity here. 
\begin{proof}[Proof of Theorem \ref{b}]
Let $\Gamma$ be the fundamental group of a closed arithmetic locally symmetric manifold $(M,g)$. And $(M',g')$ is another closed Riemannian manifold with negative curvature and homotopically equivalent to $M$.

Let $(\widetilde{M},\widetilde{g})$ and $(\widetilde{M'},\widetilde{g'})$ be the Riemannian universal cover of $(M,g)$ and $(M',g')$, respectively. We first identify their Gromov boundaries.

The lifting map of the homotopical equivalence map between $M$ and $M'$ is a quasi-isometry of $(\widetilde{M},\widetilde{g})$ and $(\widetilde{M'},\widetilde{g'})$. This map extends to a $\Gamma$-equivalent homemorphism $\Phi$ of their boundaries.

We then pull back the Busemann-cocycle $\beta' $ for $(\widetilde{M'},\widetilde{g'})$ to a cocycle on the boundary of $(\widetilde{M},\widetilde{g})$, call it $\beta$. $\beta$ is given by $\beta(\gamma, \xi)=\beta'(\gamma, \Phi \xi)$. 
Notice that the pull back does not change the marked length spectrum of cocycles, since $\Phi$ maps $\gamma^+\in \partial \widetilde{M}$ to $\gamma^+\in\partial \widetilde{M'}$ for all $\gamma\in \Gamma\setminus\{1\}$. Furthermore, pullbacks of $B$-cocycles are $B$-cocycles.

Let $\alpha$ be the Busemann-cocycle for the locally symmetric manifold $(M,g)$ with base-point $[e]$. $\alpha$ is a restriction of a cocycle $\bar{\alpha}$ of $G$. By assumption, $R_\alpha$ is a subrelation of $R_\beta$. Proposition \ref{Pro:extension} implies that $\beta$ extends to a cocycle $\bar{\beta}$ for $G$. 

Recall that the Borel $G$-cocycle on the Furstenberg boundary $X=G/P$ up to strictly equivalence is classified by $\text{Hom}(P,\bbR)=\text{Hom}(A,\bbR)$ up to equivalence \cite{zimmer2013ergodic}. In rank one case, $A=\bbR$.  Each class in $\text{Hom}(A,\bbR)$ has a $B$-cocycle representative $\lambda \bar{\alpha}$ for some $\lambda\in \bbR$.

It follows that the continuous $B$-cocycle $\bar{\beta}$ is strictly equivalent to $\lambda \bar{\alpha}$ for some $\lambda\in \bbR$. It is clear that marked length spectrum of cocycle is an invariant of strictly equivalence relation. We conclude that there exist $\lambda$ such that $\lambda \ell_g=\ell_{g'}$. In fact, by marked length spectrum rigidity for cocycles, $\beta-\lambda\alpha=d \varphi$ for some continuous map $\varphi$. 

Then the Riemannian manifold $(M',\frac{1}{\lambda}g')$ have the same marked length spectrum as $(M,g)$. By marked length spectrum rigidity for locally symmetric manifolds \cite{hamenstadt1997cocycles}, $(M,g)$ and $(M',\frac{1}{\lambda}g')$ are isometric.
\end{proof}
\section{hyperbolic surfaces without marked length pattern rigidity}\label{hypernolic surfaces}
In this section, we will show that most finite volume complete hyperbolic surfaces do not have marked length pattern rigidity.
\subsection{Fricke moduli} 
\

We first introduce a good coordinate system of the Teichmuller space for our purpose. It is called Fricke space in \cite{imayoshi2012introduction}. 

Let $S$ be a surface with genus $g$ and $n$ punctures where $2-2g-n<0$. Denote $\Gamma$ the fundamental group  of $S$, and $T_{g,n}$ the Teichmuller space of $S$. A hyperbolic structure on $S$ induce an lattice embedding of $\Gamma$ into  $G:=\PSL_2(\bbR)$. 

Let \{$\alpha_i$, $\beta_i$, $\gamma_j$\}, ($1\leq i\leq g$, $1\leq j\leq n$) be a canonical system of generators of $\Gamma$ with the fundamental relation
\[
\prod_{i=1}^g[\alpha_i,\beta_i]\prod_{j=1}^n\gamma_j=e.
\] 
Assume $\pi$ is an lattice embedding of $\Gamma$ into $G$. Up to conjugacy, we impose the normalization conditions: 
\begin{enumerate}
    \item  $\pi(\alpha_1)$ has its repelling and attractive fixed points at 0 and $\infty$, respectively,

 \item $\pi(\beta_1)$ has a fixed point at 1. 
\end{enumerate}
Then the matrix representation of  $\pi(\alpha_1)$ and $\pi(\beta_1)$ are given by 
\[
\begin{array}{ll}
\begin{pmatrix}
\lambda&0\\
0&\lambda^{-1}
\end{pmatrix},& \lambda>1,\\
\begin{pmatrix}
a_1&b_1\\
c_1&d_1
\end{pmatrix}, &a_1d_1-b_1c_1=1, \quad a_1+b_1=c_1+d_1>0
\end{array}
\]
respectively.

For $2\leq i\leq g$, $\pi(\alpha_i)$ and $\pi(\beta_i)$ are represented uniquely by the matrices 
\[
\begin{array}{lll}
\begin{pmatrix}
a_i&b_i\\
c_i&d_i
\end{pmatrix},&a_i d_i-b_i c_i=1,&c_i>0,\\
\begin{pmatrix}
a'_i&b'_i\\
c'_i&d'_i
\end{pmatrix},&a'_i d'_i-b'_i c'_i=1,&c'_i>0,
\end{array}
\] 
Similarly, for $1\leq j\leq n$, $\pi(\gamma_j)$ is written uniquely in the form
\[\begin{pmatrix}
e_j&f_j\\
g_j&h_j
\end{pmatrix}\] with $e_j h_j-f_j g_j=1$, $e_j+h_j=2$.

We define the Fricke coordinates by assign a lattice embedding to the sequence 
$(a_i, c_i, d_i, a'_i, c'_i, d'_i, e_j, g_j)$, $2\leq i\leq g$, $1\leq j\leq n$. In \cite{imayoshi2012introduction}, it is showed that the Fricke coordinates defines an embedding of the Teichmuller space $T_{g,n}$ into $\bbR^{6g-6+2n}$.

We recall the algorithm to recover $\pi$ from its Fricke coordinates. For more details, see \cite{imayoshi2012introduction}. 

It is clear that $\pi(\alpha_i)$, $\pi(\beta_i)$ and $\pi(\gamma_j)$, $2\leq i\leq g$, $1\leq j\leq n$ are uniquely determined by the Fricke coordinates. Note that all the $b_i$, $b'_i$, $2\leq i\leq g$ are rational functions of the coordinates. The same is true for $f_j$ and $h_j$, $1\leq j\leq n$.  

What remains to show is that $\pi(\alpha_1)$ and $\pi(\beta_1)$ are determined by the Fricke coordinates. 
Let \[
(\pi(\prod_{i=2}^g[\alpha_i.\beta_i]\prod_{1=1}^n\gamma_j))^{-1}=\begin{pmatrix}
a&b\\
c&d
\end{pmatrix}.
\] 
From the fundamental relation of $\Gamma$:
\[\begin{pmatrix}
\lambda&0\\
0&\lambda^{-1}
\end{pmatrix}
\begin{pmatrix}
a_1&b_1\\
c_1&d_1
\end{pmatrix}
\begin{pmatrix}
\lambda^{-1}&0\\
0&\lambda
\end{pmatrix}=\pm\begin{pmatrix}
a&b\\
c&d
\end{pmatrix}\begin{pmatrix}
a_1&b_1\\
c_1&d_1
\end{pmatrix}.\]
Replacing $\begin{pmatrix}
a&b\\
c&d
\end{pmatrix}$ by $\begin{pmatrix}
-a&-b\\
-c&-d
\end{pmatrix}$, if necessary, we have 
\begin{equation}\label{E}
\lambda^2=\frac{a-1}{1-d}.  \end{equation}
And
\[a_1=\frac{b}{1-a}c_1,\]
\[d_1=\frac{c}{1-d}b_1.\]
By $a_1+b_1=c_1+d_1$, 
\[\frac{a+b-1}{1-a}c_1=\frac{c+d-1}{1-d}b_1.\]

Recall that $a_1d_1-b_1c_1=1$, it follows
\begin{equation}\label{E1}
c^2_1[\frac{bc(a+b-1)}{(1-a)^2(c+d-1)}-\frac{c+d-1}{1-d}]=1.
\end{equation}

Equation \ref{E} and \ref{E1} show that $\lambda^2$ and $c_1^2$ are rational functions of the Fricke coordinates.

From now on, we identify Teichmuller space as a subset of $\bbR^{6g-6+2n}$ via Fricke coordinates.

\subsection{Length of geodesic and Horowitz's Theorem}
\

For any hyperbolic element $A$ in $\PSL_2(\bbR)$, let $\ell_A$ be the translation length of $A$. Then $|\text{tr}(A)|=2\cosh(\ell_A)$. Hence the marked length pattern of a Fuchsian group is determined by traces of its elements.  

Let $\pi$ be a lattice representation of the surface group $\Gamma$ into $\PSL_2(\bbR)$. Denote the Fricke coordinates of $\pi$ by $X$. Recall that $\lambda$ and $c_1$ are determined by equations \ref{E}, \ref{E1} and the normalization conditions.  

Let $s$, $t$ be two real numbers with $st\neq 0$. Then there is a family of representations $\pi'_{s,t}: F_{2g+n}\rightarrow \GL_2(\bbR),$

\[\pi_{s,t}'(\alpha_1)=\begin{pmatrix}
s&0\\
0&s^{-1}
\end{pmatrix},\]
\[\pi_{s,t}'(\beta_1)=\begin{pmatrix}
\frac{bt}{1-a}&\frac{(1-d)(a+b-1)t}{(1-a)(c+d-1)}\\
t&\frac{(a+b-1)ct}{(1-a)(c+d-1)}
\end{pmatrix},\]
\[\pi_{s,t}'(\alpha_i)=\begin{pmatrix}
a_i&\frac{a_i d_i-1}{c_i}\\
c_i&d_i
\end{pmatrix},\]
\[\pi_{s,t}'(\beta_i)=\begin{pmatrix}
a'_i&\frac{a'_id'_i-1}{c;_i}\\
c;_i&d'_i
\end{pmatrix},\]
\[\pi_{s,t}'(\gamma_j)=\begin{pmatrix}
e_j&\frac{2e_j-e^2_j-1}{g_j}\\
g_j&2-e_j
\end{pmatrix},\]
$2\leq i\leq g$, $1\leq j\leq n$ where $a$, $b$, $c$, $d$ is given by 
\[\begin{pmatrix}
a&b\\
c&d
\end{pmatrix}\prod_{i=2}^g[\pi_{s,t}'(\alpha_i),\pi_{s,t}'(\beta_i)]\prod_{j=1}^n\pi'_{s,t}(\gamma_j)=\pm\begin{pmatrix}
1&0\\
0&1
\end{pmatrix}\]
as before.
When $s=\lambda$ and $t=c_1$, $\pi'_{\lambda,c_1}$ induces the representation $\pi$.

It is clear that all the traces of $\pi_{s,t}'(w)$ where $w\in F_{2g+n}$ are rational functions of $X$, $s$ and $t$. An induction shows that for any word $w\in F_{2g+n}$, $\text{tr}(\pi'_{s,t}(w))$ has the form
\[\text{tr}(\pi'_{s,t}(w))=t^k\sum_{l=-\infty}^\infty \omega_{i}s^l \]
where $\omega_{i}$ are rational functions on $X$, $k\in \bbZ$ and all but finite many of $\omega_{i}$ are 0.

Fix $s$ and $t$. There are four different representations $\pi'_{\pm s,\pm t}$ all with the same traces up to sign. Indeed, we just replace $\pi_{s,t}'(\alpha_1)$ and $\pi'_{s,t}(\beta_1)$ by $\pm\pi'_{s,t}(\alpha_1)$ and $\pm\pi'_{s,t}(\beta_1)$, respectively. It follows that $(\text{tr}(\pi'_{\pm s,\pm t}(w)))^2$ is independent of the choice of these four representations. Hence $(\text{tr}(\pi'_{s,t}(w)))^2$ has the form 
\[(\text{tr}(\pi_{s,t}'(w)))^2=t^{2k}\sum_{l=-\infty}^\infty \Omega_{i}s^{2l} \]
where $\Omega_{i}$ are rational functions on $X$, and all but finite many of $\Omega_{i}$ are 0.

Let $s=\lambda$, $t=c_1$. $\pi'_{\lambda,c_1}$ is indeed a representation of $F_{2g+n}$ to  $\SL_2(\bbR)$. By \ref{E} and \ref{E1}, $(\text{tr}(\pi_{\lambda, c_1}'(w)))^2$ is a rational function on $X$ for all $w\in F_{2g+n}$. Moreover,
 \[
 \pi_{\lambda, c_1}'(\prod_{i=1}^g[\alpha_i,\beta_i]\prod_{j=1}^n\gamma_j)=\pm \begin{pmatrix}
1&0\\
0&1
\end{pmatrix}.
\]
Hence $\pi_{\lambda, c_1}'$ is a lifting of the representation $\pi$. Let $p$ be the natural projection of $F_{2g+n}$ to $\Gamma$. The traces of $\pi$ and $\pi'_{\lambda,c_1}$ are the same up to sign, i.e., for any $w\in F_{2g+n},$
\begin{equation}\label{eq}
    {\rm{tr}(\pi_{\lambda,c_1}'(w))}=\pm {\rm{tr}}(\pi(p(w))).
\end{equation}

We conclude that $({\rm{tr}}(\pi(\gamma)))^2$ are rational functions on $X$ for all $\gamma\in \Gamma$. Since Teichmuller space is an open subset of $\bbR^{6g-6+2n}$ in the Hausdorff topology, it is Zariski dense. There are unique extensions of $({\rm{tr}}(\pi(\gamma)))^2$ as rational functions on $\bbR^{6g-6+2n}$. We call the function $({\rm{tr}}(\pi(\gamma)))^2$ the rational function of $\gamma$, and denote it by $Q_\gamma$.

There is a stronger version from a theorem of Horowitz \cite{horowitz1972characters}. 
\begin{theorem}
Let $F=\langle s_1,s_2,\cdots, s_m\rangle$ be a free group on $m$ generators. For any word $w\in F$. There is a polynomial $P_w$ depends only on $w$ with integer coefficients in $2^m-1$ characters such that for any representation $\phi: F\rightarrow \SL_2(\bbR)$, 
\[{\rm{tr}}(\phi(w))=P_w(t_1,t_2,\cdots, t_{12},\cdots,t_{12\cdots m}),\]
where $t_{i_1i_2\cdots i_v}={\rm{tr}}(\phi(s_{i_1}s_{i_2}\cdots s_{i_v}))$, $1\leq i_1<i_2<\cdots<i_v\leq m$.
\end{theorem}
 
Define a relation $R_{\rm min}$ on $\Gamma$ as follows:
$\gamma R_{\rm min}\eta$ when $Q_\gamma=Q_\eta$. It is clear now that for any hyperbolic metric $g$ on $S$, $R_{\rm min}$ is a sub-relation of $R_g$. 

 For any $\gamma$ and $\eta\in \Gamma$, the rational equation $Q_\gamma=Q_\eta$ defines a algebraic subset of $\bbR^{6g-6+2n}$. The intersection of this set and the Teichmuller space is either everything or a subset of positive codimension. Since there are just countably many pair of elements, we get countably many algebraic subsets of positive codimension. By dimension reasons or refer to Lebesgue measure, we know that the union of all these countable subsets of positive codimension $T_{\mathrm{singular}}$ is a proper subset of the Teichmuller space. 

Choose two points in the complimentary of $T_{\mathrm{singular}}$. We have two hyperbolic metrics $g_1$ and $g_2$, with $R_{\rm min}=R_{g_1}=R_{g_2}$. Hence $g_1$ do not share marked length pattern rigidity. 
\appendix
\section{General geometric boundaries}\label{appendix}
We generalize our result to general acylindrical groups in this appendix. The proof are almost the same as in the paper. We just point out some necessary changes and restate the theorems in this case. 

General acylindrical groups and their Possion-Furstenberg boundaries are not necessarily geometric boundary, since there are elliptic and parabolic elements. To deal with it, we introduce the following notation.

Let $X$ be a $\Gamma$-space and $\gamma\in \Gamma$ acts on $X$ hyperbolicly. For any $\eta\in \Gamma$ such that $\eta \gamma^+\neq \gamma^-$, let 
\[
A_{\eta,\gamma}=\setdef{n}{\eta\gamma^n \ is\ hyperbolic.}=\setdef{n}{\gamma^n\eta \ is\ hyperbolic.}.
\]
The last equation is true because $\eta\gamma^n=\eta(\gamma^n\eta)\eta^{-1}.$ Notice that $A_{\eta,\gamma}=A_{\eta^{-1},\gamma^{-1}}$.
\begin{defn}
Let $\Gamma$ be a topological group. A nontrivial compact Hausdorff $\Gamma$-space $X$ is call a \textbf{general geometric boundary} if
\begin{enumerate}
    \item[(1)]  $\Gamma$ acts on $X$ minimally,
    \item[(ii)] There are hyperbolic elements and $\sup A_{\eta,\gamma}=+\infty$ for all hyperbolic element $\gamma$ when $\eta \gamma^+\neq \gamma^-$,  
    \item[(3)] There are $\Gamma$-quasi invariant measures $\mu$, $\mu'$ on $X$ such that $\Gamma$is $\mu\times \mu'$-ergodic. 
\end{enumerate}
If in addition, $\mu=\mu'$, we call $X$ a symmetric general geometric boundary. 
\end{defn}

We just replace (2) in the definition of geometric boundary by $(ii)$. 

First, we show that general acylindrical groups and their Possion-Furstenberg boundaries are general geometric boundaries.

Recall that a group $G$ is called acylindrically hyperbolic if $G$ admits a non-elementary acylindrical isometric action on a geodesic Gromov-hyperbolic space $M$. It was showed by Maher and Tiozzo in \cite{maher2018random} that the Furstenberg-Poisson boundary of a spread-out generating measure on $G$ is the same as the limit set of $G$ in $\partial M$ with the hitting measure. By the work of Bader and Furman \cite{bader2014boundaries}, (3) follows. (1) is true since $G$ acts on its limit set minimally. 

For $(ii)$, recall that there is visual metric $d$ on the limit set $X$. For hyperbolic element $\gamma$, and any compact subset $K\subset X-\{\gamma^-\}$, there exist $L$, $\kappa>0$ such that $d(\gamma^n x,\gamma^+)\leq L\exp(-n\kappa)$ for all $x\in K$. Let $\eta\in \Gamma$. $\eta$ acts on $(X,d)$ by Lipschitz homemorphism. Take an open neighbourhood $U$ of $\eta\gamma^+$ with $\gamma^-\notin U$. There is a $N$ such that for all $n\geq N$, $\gamma^n U\subset \eta^{-1}U$. Increasing $n$ if necessary, $\eta\gamma^n|_U$ is a contraction. Hence there is a contracting fixed point of $\eta\gamma^n$ in $U$. By classifying of isometries of hyperbolic spaces, $\eta\gamma^n$ is hyperbolic for all $n$ big enough. 

Second, all lemmas in section \ref{4} still hold. The proofs are the same, we restate them here. 

\begin{lemma}
General geometric boundaries have infinite many points.
\end{lemma}

\begin{lemma}
Let $\Gamma$-space $X$ be a general geometric boundary. Let $\gamma$ be a hyperbolic element. There exist $\theta\in \Gamma$ such that $\theta \gamma^+\neq \gamma^-$, $\theta \gamma^-\neq \gamma^-$.
\end{lemma}
\begin{lemma}
Conjugations of a hyperbolic element are hyperbolic. And $$(\theta\gamma\theta^{-1})^\pm=\theta\gamma^\pm$$
when $\gamma$ is hyperbolic.
\end{lemma}
\begin{lemma}
Let $\gamma$ and $\eta\in \Gamma$ with $\gamma$ hyperbolic. If $\eta \gamma^+\neq \gamma^-$, then 
$$\lim_{n\in A_{\eta,\gamma},n\rightarrow \infty}(\eta\gamma^n)^+=\eta\gamma^+,$$
$$\lim_{n\in A_{\eta,\gamma},n\rightarrow \infty}(\eta\gamma^n)^-=\gamma^-.$$
\end{lemma}

Define the \emph{marked length spectrum function} 
for general cocycle $\beta$ on a general geometric boundary $X$ 
by setting
\begin{equation}
    \ell_\beta(\langle \gamma\rangle)=\beta(\gamma,\gamma^+)
\end{equation}
for all hyperbolic elements $\gamma\in \Gamma\setminus\{1\}$.

The same proof of Theorem \ref{a} using these lemmas gives
\begin{theorem}\label{A.6}
Let $\Gamma$-space $X$ be a general geometric boundary and $\alpha, \beta: \Gamma\times X\rightarrow \bbR$
two $B$-cocycles. If $\ell_\alpha=\ell_\beta$, then $\alpha-\beta=d \varphi$ for some continuous function $\varphi$.
In other words, $[\alpha]=[\beta]$ in $H^1_c(\Gamma,X,\bbR)$.
\end{theorem}

We generalise the definition of Marked length spectrum function to finite volume negatively curved manifolds.

Let $(M,g)$ be a complete finite volume manifold with negative curvature. The fundamental group $\Gamma=\pi_1(M)$ contains parabolic elements. Let $\gamma \in \Gamma$ be parabolic, there is arbitrary short closed geodesic represents $\gamma$. We take the convention that $\ell_g(\langle\gamma
\rangle)=0$ for parabolic classes. It was called minimal marked length spectrum in some papers. It is the infimum of the length of all closed geodesics in class $\langle\gamma\rangle$. We define the marked length pattern by the same way as before:
\[
    R_g=\setdef{(c_1,c_2)\in C_\Gamma\times C_\Gamma}{\ell_g(c_1)=\ell_g(c_2)}.
\]

As a application of Theorem \ref{A.6}. By the same construction in this paper, we are able to show the following two Theorems.
\begin{theorem}\label{A.7}
Let $(M,g_1)$ and $(M,g_2)$ be two arbitrary finite volume complete closed strictly negatively curved  Riemannian metrics on a manifold $M$ with fundamental group $\Gamma$. Let $H$ be a subgroup of $\Gamma$ such that the limit set of $H$ is all of $\partial \widetilde{M}$. 

Then, $\ell_{g_1}=\ell_{g_2}$ on classes from $H$ 
only if $\ell_{g_1}=\ell_{g_2}$ on all of $\Gamma$. 
\end{theorem}

\begin{theorem}\label{A.8}
    Let $(M,g_0)$ be a finite volume arithmetic locally symmetric manifold of rank 1, and let $g$ 
    be an arbitrary strictly negatively curved complete finite volume Riemannian metric on $M$.
    Then $R_{g_0}\subset R_g$ only if $\ell_{g_0}=\lambda\ell_{g}$ for some $\lambda>0$.
\end{theorem}

\begin{remark}
In \cite{cao1995rigidity}, Cao showed that if
two orientable, uniform visibility surfaces of finite area and bounded non-positive
curvature have the same marked length spectrum, then they must be isometric. Hence we can strength Theorem \ref{A.7} and \ref{A.8} in dimension 2.

In higher dimension, we do not know the marked length spectrum rigidity.  Peigné and Sambusetti \cite{peigne2019entropy} showed the following:

\emph{Let $M$ be a finite volume $n$-manifold with pinched, negative curvature
$-b^2\leq \kappa\leq -1$ which is homotopy equivalent to a locally symmetric manifold $M_0$ with curvature normalized between -4 and -1. If $M$ and $M_0$ has same marked length spectrum, then they are isometric.}
\end{remark}

\end{document}